\documentclass{article}
\usepackage{amsfonts,amsmath,amstext,amsbsy,euscript,amssymb, graphicx}
\begin{document}
\title{Catalan rook and queen paths with boundary}
%
\author{Joseph P.S. Kung \thanks{Department of Mathematics, University of North Texas, Denton TX 76203, USA. \texttt{kung@unt.edu}}, Anna de Mier\thanks{Departament de Matem\`atica Aplicada II, Universitat Polit\`ecnica de Catalunya, 08034 Barcelona, Spain. \texttt{anna.de.mier@upc.edu}}}

\maketitle
\abstract{
This is the second version of the paper.  There are several additional results.  This version has been submitted for publication under the title ``Catalan lattice paths with rook, queen and spider steps.''  
 
\vskip 0.2in
A lattice path is a path on lattice points (points with integer coordinates) in the plane in which any step increases the $x$- or $y$-coordinate, or both. A rook step is a proper horizontal step east or vertical step north.  A bishop step is a proper diagonal step of slope $1$ (to the northeast).   A spider step is a proper step of finite slope greater than $1$ (in a direction between north and northeast).  
A lattice path is Catalan if it starts at the origin and stays strictly to the left of 
the line $ y = x - 1.$  
We give abstract formulas for the ordinary generating function of the number of lattice paths with a given right boundary and steps satisfying a natural slope condition.  Explicit formulas are derived for generating functions of the number of 
Catalan paths in which all rook steps and some (or all) bishop or spider steps are allowed finishing at $(n,n).$  
These generating functions are algebraic; indeed, many satisfy quadratic equations.  }

\newpage
\vskip 0.4in\noindent
{\bf  1.  Paths and boundaries}  

Let $S$ be a subset of $\mathbb{N} \times \mathbb{N},$ where $\mathbb{N}$ is the set of non-negative integers.  
An {\sl $S$-path} from the point $(n_0,m_0)$ to the point $(n,m)$  is a sequence of pairs (called {\sl steps}) $(a_1,b_1), (a_2,b_2), \ldots, (a_k,b_k)$ from $S$
such that 
$$
(n_0 +a_1 + a_2 + \cdots + a_k, m_0 + b_1 + b_2 + \cdots + b_k)  = (n,m).  
$$
The steps $(a,0)$ and $(0,a),$ $a > 0,$ are called {\sl rook} steps, the steps $(a,a), \, a > 0,$ are called 
{\sl bishop} steps, and the steps $(a,b),\, b > a > 0$ are called {\sl spider} steps.  
Let $a_{n,m}$ be the number of $S$-paths from the origin $(0,0)$ to the point $(n,m).$   
Then $a_{n,m}$ is finite if and only if $(0,0) \not\in S.$  If $(0,0) \notin S,$ then a simple argument yields the following bivariate generating function 
$$
A(x,y) = \sum_{n,m = 0}^{\infty} a_{n,m} x^ny^m = \frac {1}{1 - \sum_{(a,b) \in S} x^a y^b}.  
$$ 
In the rest of this paper, we shall always assume that $(0,0) \notin S.$    
Examples of $S$-paths are {\sl classical  lattice paths} ($S = \{(1,0),(0,1)\}$), {\sl rook paths} ($S = \{(a,0),(0,a): a > 0 \}$), and {\sl queen paths} ($S = \{(a,0),(0,a): a > 0\} \cup \{(a,a): a>0\}$).  

We can represent an  $S$-path $(a_1,b_1), (a_2,b_2), \ldots, (a_k,b_k)$   
from $(n_0,m_0)$ to $(n,m)$ geometrically by the union of the line segments
\begin{eqnarray*} 
&&  \overline{(n_0,m_0),(n_0+a_1,m_0+b_1)}, \overline{(n_0+a_1,m_0+b_1), (n_0+a_1+a_2,m_0+b_1+b_2)},\ldots, 
\\   
&&  \qquad\qquad\qquad\qquad   
\overline{(n_0+a_1+a_2 + \cdots + a_{k-1},m_0+b_1+b_2 + \cdots + b_{k-1}),(n,m)},  
\end{eqnarray*}
so that an $S$-path is a piecewise linear path.   
We call the points 
$$ 
(n_0,m_0),(n_0+a_1,m_0+b_1), (n_0+a_1+a_2,m_0+b_1+b_2),\ldots, $$ 
$$(n_0+a_1+a_2 + \cdots + a_{k-1},m_0+b_1+b_2 + \cdots + b_{k-1}),(n,m),  
$$
the {\sl nodes} of the path.  

Let $\underline{s}$ be a non-decreasing sequence $(s_i)_{0 \le i < \infty}$ of positive integers.  An $S$-path starting from $(0,0)$ {\sl has (right) boundary} $\underline{s}$ if as a piecewise linear path, 
it lies strictly to the left 
of the  union of the line segments  
$$
\overline{(s_0,0), (s_1,1)},  \overline{(s_1,1), (s_2,2)},  \overline{(s_2,2),
  (s_3,3)},\ldots \,,  $$
or, equivalently, every node of the $S$-path of the form $(x,i)$ satisfies $x<s_i$.

An $S$-path is {\sl Catalan} if it starts at the origin $(0,0)$ and has right boundary $(1,2,3,\ldots).$  
When a Catalan $\{(0,1),(1,0)\}$-path is rotated clockwise by $45$ degrees, one obtains a Dyck path.  Thus, Catalan $S$-paths may be considered as Dyck paths with more types of steps.      

A set $S$ satisfies the {\sl slope condition for the boundary $\underline{s}$} if 
for every pair $(a,b)$ with $b \neq 0$ in $S,$ and every index $i,$ 
$$
s_{i+1} > s_i - 1 +a/b, s_{i+2} > s_i - 1 + 2a/b, \ldots, s_{i+b} > s_i - 1 +a, 
$$
or, equivalently, an $S$-path can touch or cross the boundary $\underline{s}$ only at a horizontal step.    
Classical lattice paths, rook paths, and queen paths all satisfy the slope condition relative to 
the Catalan boundary $(i+1).$

\vskip 0.4in\noindent 
{\bf  2.  A combinatorial bijection} 

\vskip 0.2in\noindent 

Let $\underline{s} = (s_i),$ be a non-decreasing sequence of positive integers, $ \mathrm{AP}(n,m)$ be the set of all $S$-paths from $(0,0)$ to $(n,m),$ $a_{n,m} = |\mathrm{AP}(n,m)|,$  
$\mathrm{Path}(n)$ be the set of $S$-paths from $(0,0)$ to $(s_n-1,n)$ with right boundary $\underline{s},$  
$p_n = |\mathrm{Path}(n)|,$ and $\mathrm{LHP}(n)$ be the subset of $S$-paths in $\mathrm{Path}(n)$ such that 
the last step is a horizontal step.

We begin with a combinatorial bijection.  The idea behind the bijection was
used in an earlier paper \cite{JSIP}. 

\vskip 0.2in\noindent
{\bf Lemma 2.1.}  Suppose that $S$ satisfies the slope condition for the
boundary $\underline{s}$ and that    $S$ contains all proper horizontal steps $(a,0), \, a > 0.$  Then there is a bijection between 

$$ \mathrm{AP}(s_n-1, n)
$$ and the following union of three disjoint subsets:

\begin{align*} 
& \bigcup_{m=0}^{n-1}
\left[\mathrm{Path}(m) \times  [\cup_{j=0}^{s_n-s_m-1} \mathrm{AP}(j,n-m)] 
\right]
\\& 
\cup \, 
\bigcup_{m=0}^{n-1}
\left[\mathrm{LHP}(m) \times  [\cup_{j=0}^{s_n-s_m-1} \mathrm{AP}(j,n-m)] 
\right]
\\
& \cup \,
\mathrm{Path}(n).  
\end{align*}

\vskip 0.1in\noindent

{\bf Proof.}  
Paths in $\mathrm{AP}(s_n-1,n)$ either have boundary $\underline{s}$ or not.  Those having boundary $\underline{s}$ are mapped onto themselves in the third subset $\mathrm{Path}(n).$  

Next consider $S$-paths from $(0,0)$ to $(s_n-1,n)$ which fail to have 
boundary $\underline{s}.$  We will decompose these paths into an initial subpath and a final subpath. Let $(d,m)$ be the first node on or to the right of the boundary.  Then $0\leq m\leq n-1$ and $s_m \leq d \leq s_n-1.$   From the point $(d,m)$, the path goes as it wishes to the point $(s_n-1,n)$, so that the set of final subpaths from $(d,m)$ to $(s_n-1,n)$ is in bijection, by a translation, with the set $\mathrm{AP}(s_n-d-1,n-m).$  

There are more possibilities for the initial subpath from $(0,0)$ to $(d,m).$  By the slope condition, the step to $(d,m)$ is a horizontal step.  As horizontal steps can be of any nonzero length, the node before $(d,m)$ was $(c,m)$ for some $c\leq s_m-1.$ There are two cases, according to whether this inequality is an equality or not (the latter being only possible if $s_m>1$).

\vskip 0.1in

{\it Case 1, $c = s_m-1.$}  The initial subpath consists of a path in $\mathrm{Path}(m),$ followed by the step $(d-s_m+1,0),$ taking us to $(d,m).$  In this case, the path is decomposed into an initial subpath in $\mathrm{Path}(m)$ followed by a horizontal step of size $d-s_m +1$ and a final subpath which is in $\mathrm{AP}(s_n-d-1,n-m)$ after a translation.  
\vskip 0.1in

{\it Case 2,  $c < s_m-1.$}  The last step in the initial subpath is a ``long'' horizontal step taking us from $(c,m)$ through the boundary to $(d,m).$   We break up this long step into two shorter steps at $(s_m-1,m).$  Specifically, we take the path from $(0,0)$ to $(c,m),$ and
add a {\sl horizontal} step $(s_m-c-1,0),$ so that we now end at $(s_m-1,m).$      We then continue with a step $(d-s_m+1,0)$ which brings us through the boundary to $(d,m).$   In this case, the path is decomposed into an initial subpath in $\mathrm{LHP}(m)$ and a final subpath which is in
$\mathrm{AP}(s_n-d-1,n-m)$ after a translation.

In both cases, the decomposition is reversible and hence, bijective.   This completes the proof of the lemma.


\vskip 0.2in
Lemma 2.1 implies immediately that  
\begin{align*}
a_{s_n-1,n} &= p_n+\sum_{m=0}^{n-1}   p_m(a_{0,n-m}+a_{1,n-m}+\cdots+a_{s_n-s_m-1,n-m})\\
 &+ \sum_{m=0}^{n-1}   |\mathrm{LHP}(m)| (a_{0,n-m}+a_{1,n-m}+\cdots+a_{s_n-s_m-1,n-m}).
\end{align*}

From now on we restrict to the case $s_0=1$. Then $p_0=1$ and
$|\mathrm{LHP}(0)|=0$. For $m\geq 1$, we write $|\mathrm{LHP}(m)|=p_m-d_m$,
where  $d_m$ is the number of paths in $\mathrm{Path}(m)$ that end in a
non-horizontal step. 
 By doing minor manipulations and setting $d_0=0$, we obtain the following corollary.

\vskip 0.2in\noindent 

{\bf Corollary 2.2.}  Under the hypotheses in Lemma 2.1, if $s_0=1$ then  
\begin{align*}
& a_{0,n}+a_{1,n}+\cdots+a_{s_n-2,n} + a_{s_n-1,n}  \\
& = p_n+ \sum_{m=0}^{n-1}   2p_m(a_{0,n-m}+a_{1,n-m}+\cdots a_{s_n-s_m-1,n-m}) \\
 & \quad - \sum_{m=0}^{n-1}  d_m (a_{0,n-m}+a_{1,n-m}+\cdots+a_{s_n-s_m-1,n-m}).
\end{align*}

In the Catalan case, where the boundary is $( i+1),$ the recursions in Corollary 2.2 can be combined to obtain a formula for the ordinary generating function of $p_n.$  
Let $h$ be an integer and 
\begin{align*}
P(t) &= \sum_{n=0}^{\infty} p_n t^n, \\
Q(t) &= \sum_{n=0}^{\infty} d_n t^n, \\  
D_h(t)  &=\sum_{n=h}^{\infty} a_{n-h,n} t^n.
\end{align*}
 In general, $D_h(t)$ the generating function of $S$-paths in a rectangle where the height is $h$ units ``larger'' than the width.   In particular, $D_0(t)$ is the diagonal of the bivariate generating function $A(x,y)$ for $S$-paths. 

To calculate the series $Q(t)$, we observe that the slope condition for the Catalan boundary implies that if a path in $\mathrm{Path}(m)$ ends with a non-horizontal step, then that step must be a bishop step $(a,a)$ for some $a > 0.$  Removing that bishop step yields a path in $\mathrm{Path}(m-a).$  Conversely, if 
$(b,b) \in S,$ every path in $\mathrm{Path}(m-b)$ extends to a path in $\mathrm{Path}(m)$ with last step $(b,b).$ Hence, 
$$
d_m = \sum_{a \in I} p_{m-a}, 
$$
where $I =  \{a: (a,a) \in S\}.$  In terms of generating functions, this says that $Q(t) = P(t)T(t),$ where  
$$
T(t) = \sum_{a\in I} t^a.    
$$
For rook paths, $T(t) = 0,$ and for queen paths, $T(t) = \frac {t}{1-t}.$

Now take the recursions in Corollary 2.2 and multiply the $n$-th recursion by $t^n,$ obtaining 
\begin{align*}
a_{0,0} &=  p_0
\\
(a_{0,1} + a_{1,1})t  &=  
p_1 t+  2p_0 a_{0,1}t - d_0a_{0,1}t 
\\ 
(a_{0,2} + a_{1,2} + a_{2,2}) t^2 &=  p_2 t^2 + 2 p_1 a_{0,1}t^2 + 2p_0 (a_{0,2}+a_{1,2})t^2  -  d_1 a_{0,1}t^2 
- d_0 (a_{0,2}+a_{1,2})t^2
\\
(a_{0,3} + a_{1,3} + a_{2,3} + a_{3,3}) t^3 &=  
p_3 t^3 + 2p_2 a_{0,1} t^3 + 2 p_1( a_{0,2} + a_{1,2})t^3 +  2p_0 (a_{0,3}+a_{1,3} + a_{2,3})t^3
\\
&\quad  - d_2 a_{0,1} t^3 - d_1( a_{0,2} + a_{1,2})t^3 -d_0 (a_{0,3}+a_{1,3} + a_{2,3})t^3, 
\\
& \vdots   
\end{align*}
with the $n$-th equation given by  
\begin{align*} 
(a_{0,n} +a_{1,n}  + \cdots + a_{n,n}) t^n &=  
\\
 p_n  t^n + 2p_{n-1}a_{0,1}t^n  &+   2p_{n-2}(a_{0,2}+a_{1,2})t^n  +
\cdots  + 2p_0 (a_{0,n}+a_{1,n}+ \cdots + a_{n-1,n})t^n
\\
- d_{n-1}a_{0,1}t^n  &-   d_{n-2}(a_{0,2}+a_{1,2})t^n  -
\cdots  -d_0 (a_{0,n}+a_{1,n}+ \cdots + a_{n-1,n})t^n.
 \end{align*} 
Summing down subdiagonals on the left hand side, we obtain 
$$
D_0(t) + D_1(t) + D_2(t) + \cdots 
$$
For the right hand side,  observe that 
\begin{align*}
P(t) D_1(t)   
%
%
&=  p_0a_{0,1}t + (p_0 a_{1,2} + p_1 a_{0,1})t^2 + (p_0a_{2,3} + p_1 a_{1,2} + p_2a_{0,1})t^3 + \cdots \,,
\\ 
P(t) D_2(t)   
&=  p_0a_{0,2}t^2 + (p_0 a_{1,3} + p_1 a_{0,2})t^3 + (p_0a_{2,4} + p_1 a_{1,3} + p_2a_{0,2})t^4 + \cdots \,, 
\end{align*} 
and so on.  Analogous formulas hold for $Q(t):$     
$$ Q(t) D_1(t)   
=  d_0a_{0,1}t + (d_0 a_{1,2} + d_1 a_{0,1})t^2 + (d_0a_{2,3} + d_1 a_{1,2} + d_2a_{0,1})t^3 + \cdots \,,
$$
and so on.  Summing along ``subdiagonals'' on the right hand side, we conclude that 
\begin{align*}
& D_0(t)+D_1(t)+D_2(t)+\cdots \\
&= 
P(t)+ 2P(t)(D_1(t)+D_2(t)+\cdots) - Q(t)(D_1(t)+D_2(t)+\cdots). 
\end{align*}
Using the fact that $Q(t) = P(t)T(t),$ we obtain the following theorem. 

\vskip 0.2in\noindent  
{\bf Theorem 2.3.}  If $S$ contains all proper horizontal steps $(a,0), a > 0,$ and satisfies the slope condition for the boundary $(i+1)$,  then 
$$
P(t)=\frac{D_0(t)+D_1(t)+D_2(t)+\cdots}{1+(2-T(t))(D_1(t)+D_2(t)+\cdots)}.   
$$

\vskip 0.2in
Our method can be adapted to the following theorem.    

\vskip 0.2in\noindent  
{\bf Theorem 2.4.}  
If $S$ contains $(1,0)$ and no other horizontal step and satisfies the slope condition for the boundary $(i+1)$,  then 
$$
P(t)=\frac{D_0(t)}{1 + D_1(t)}.   
$$ 

\vskip 0.2in 
Theorems 2.3 and 2.4 cover the extremes, where $S$ contains exactly the unit horizontal step or all horizontal steps to the right.  Let $M$ be a positive integer and, for $a<M$, let $\mathrm{LHP}_a(n)$ be the subset of those $S$-paths in $\mathrm{Path}(n)$ that end in a horizontal step $(a,0)$.     

\vskip 0.2in\noindent
{\bf Lemma 2.5.}  Let $M$ be an integer greater than $1.$  Suppose that
$S$ satisfies the slope condition for the boundary $\underline{s}$ and that     $S$ contains exactly the horizontal steps $(a,0), \, 0 < a < M.$  Then there is a bijection between $ \mathrm{AP}(s_n-1, n) $
and the following union of three disjoint subsets:
\begin{align*} 
& \bigcup_{m=0}^{n-1} 
\left[\mathrm{Path}(m) \times  [\cup_{j=s_n-s_m-M}^{s_n-s_m-1} \mathrm{AP}(j,n-m)] 
\right]
\\ & \qquad
\cup \, 
\bigcup_{m=0}^{n-1}
\left[\bigcup_{a=1}^{M-1}\mathrm{LHP}_a(m) \times  [\cup_{j=s_n-s_m-M+a}^{s_n-s_m-1} \mathrm{AP}(j,n-m)] \right]
 \cup \, \mathrm{Path}(n).  
\end{align*}

\vskip 0.1in\noindent 

The proof of Lemma 2.5 is similar to the proof of Lemma 2.1, but we need to understand how the restriction on the length of the horizontal steps affects the bijection. In the bijection, each path in $AP(s_n-1,n)$ is decomposed into an initial subpath ending at $(s_m-1,m)$ for some $m$, a horizontal step, and a final subpath. In Case 1, the only effect is that the possible starting points for the final subpath are restricted. In Case 2,   two proper horizontal steps are joined into one ``long'' horizontal step.  This can only be done if the sum of the length of the last step in the initial subpath and the length of the adjoined horizontal step does not exceed $M.$ This imposes restrictions on the length of the last step as well as  the starting point of the final subpath.

The combinatorics of $\mathrm{LHP}_a(m)$ is more complicated and there seems to be no direct analog of Theorem~2.3.

\vskip 0.4in \noindent 
{\bf 3.  Catalan rook and queen paths}   
\vskip 0.2in\noindent  

Let $\mathbb{P}= \{1,2,\ldots\},$ the set of positive integers, and $I \subseteq \mathbb{P}.$  
An {\sl  $I$-queen path} is an $S$-path where $S = \{(a,0),(0,a): 0 < a < \infty\} \cup \{(a,a): a \in I\}.$   
For example, a rook path is an $\emptyset$-queen path and a queen path is a $\mathbb{P}$-queen path. 
In this section, we show how the formula in Theorem 2.3 can be made explicit for Catalan $I$-queen paths.  

As in Section 2, let $   T(t) = \sum_{a \in I}  t^a. $
Then the bivariate generating function $A(x,y)$ for the number 
$a_{n,m}$ 
of $I$-queen paths from $(0,0)$ to $(n,m)$ is given by  
$$
A(x,y)=  \left( 1-\frac{x}{1-x}-\frac{y}{1-y} - T(xy) \right)^{-1}. 
$$
The series $D_h(t)$ can be calculated from $A(x,y)$ using a standard method (see \cite{Furst, Gessel} or \cite{EC2}, Section 6.3).
This method is based on the observation that   
$D_h(t)$ is the coefficient of $s^h$ when $A(s^{-1},st)$ is expanded as a doubly infinite series in $s.$ 
Let  
$$
B(s,t)=A(s^{-1},st).
$$  
Then 
$$
B(s,t) = 
\left( C + \frac {1}{1-s} - \frac {st}{1-st} \right)^{-1}, 
$$
where $ C = 1 - T(t). $ 
Expanding, we have 
$$
B(s,t) = \frac {1}{1+C} + \frac {1}{1+C} \left[
\frac{ \displaystyle{ s \left(\frac {t-1}{1+C}\right) }  }
{\displaystyle {1 - \left( \frac{C+Ct + 2t}{1+C}\right)s + ts^2 }} 
\right].
$$
Let $\alpha$ and $\beta$ be the roots of the denominator in $B(s,t),$ with 
\begin{eqnarray*} 
\alpha &=&\frac 
{ C + (C + 2) t +
\sqrt{C^2 - (2C^2 + 4C + 4)t + (C+2)^2 t^2}}
{2(1 + C)t} 
\\
&=&   \frac {C}{1+C} t^{-1} - \frac {1}{C(1+C)} -  \frac {1+C}{C^3}t +  \cdots  
\end{eqnarray*}
and 
\begin{eqnarray*} 
\beta&=& \frac 
{ C + (C + 2) t -
\sqrt{C^2 - (2C^2 + 4C + 4)t + (C+2)^2 t^2}}
{2(1 + C)t} 
\\
&=&     
\frac {1+C}{C}  + \frac {1+C}{C^3}t + \cdots \,.
\end{eqnarray*}
Next, we expand $B(s,t)$ into a partial fraction, obtaining  
$$
B(s,t) = \frac {1}{1+C} + \frac {1}{(1+C)}\left[\frac {A}{t(s - \alpha)} + \frac {B}{t(s - \beta)}  \right],  
$$
where 
$$
A = \frac  {\alpha}{(1+C)(\alpha - \beta)} (t-1), \quad B = - \frac  {\beta}{(1+C)(\alpha - \beta)} (t-1).
$$  
Since $\alpha$ is a proper Laurent series (that is, it has negative powers of $t$), $\beta$ is a power series (that is, it has no negative powers of $t$), and $B(s,t)$ is a doubly infinite series in $s$, but not in $t$, the correct expression is found by writing 
$$
B(s,t)= \frac{1}{1 + C} -\frac{A}{(1+C) t\alpha (1-s\alpha^{-1})}  + \frac{B}{(1+C)ts(1-s^{-1}\beta)}.
$$
Expanding, we have   
$$   
B(s,t)  = \frac{1}{1 + C} -\frac{A}{(1+C) t\alpha}\left(1 + \frac {s}{\alpha}+ \frac {s^2}{\alpha^2} +\cdots \right)
+ \frac{B}{(1+C)t}\left(\frac {1}{s} + \frac {\beta}{s^2} + \frac {\beta^2}{s^3} + \cdots  \right).
$$
From this, we see easily that   
$$
D_0(t) =  \frac{1}{1+C}-\frac{A}{(1+C)t\alpha}
$$
and for $h>0,$ 
\begin{eqnarray*}
D_h(t)&=&  -\frac{A}{(1+C)t \alpha^{h+1}},  
\end{eqnarray*}
as well as 
\begin{eqnarray*} 
\sum_{h=0}^{\infty} D_h(t) &=& \frac{1}{1+C}+\frac{A}{(1+C)t(1-\alpha)}, 
\\
\sum_{h=1}^{\infty}  D_h(t)  &=& \frac{A}{(1+C)t\alpha(1-\alpha)}. 
\end{eqnarray*} 
Using Theorem 2.3, we conclude that 
$$ 
P(t) \quad=\quad 
\frac{\displaystyle{\frac{1}{1+C}+\frac{A}{(1+C)t(1-\alpha)}}}{\displaystyle{1+ \frac{(2-T)A}{(1+C)t\alpha(1-\alpha)}}}
\quad = \quad 
\frac{\displaystyle{1+\frac{A}{t(1-\alpha)}}}{\displaystyle{1+C+\frac{(2-T)A}{t\alpha(1-\alpha)}}}. 
$$

We now have the tools to derive several explicit formulas in a uniform way.  We begin with Catalan rook paths.  In this 
case, $T(t) = 0, C = 1,$ 
$$
\alpha=\frac{1+3t+\sqrt{1-10t+9t^2}}{4t}, \quad \beta=\frac{1+3t-\sqrt{1-10t+9t^2}}{4t},   
$$ 
and  
$$
A=\frac{1-10t+9t^2+(1+3t)\sqrt{1-10t+9t^2}}{4(9t-1)}.
$$
After some algebra, we obtain the following generating function.  

\vskip 0.2in\noindent  
{\bf Theorem 3.1} (Catalan rook paths). 
$$
P_{\mathrm{rook}}(t)=\frac{1+3t-\sqrt{1-10t+9t^2}}{8t}=\frac{\beta}{2}.
$$

\vskip 0.2in\noindent  
Expanding this formula gives 
$$
P_{\mathrm{rook}}(t) = 1 + t + 5t^2 + 29t^3 + 185t^4 + 1257t^5 + \cdots.
$$ 
Using standard techniques (see, for example, \cite[Chapter VI]{acFS}) we obtain the following asymptotic formula. 

\vskip 0.2in\noindent  
{\bf Corollary 3.2.} 
$$
p_{\mathrm{rook},n} \sim\frac{3\sqrt{2}}{8} \cdot \frac{9^n}{\sqrt{\pi n^3}}.
$$

\vskip 0.2in\noindent  
Note that $P_{\mathrm{rook}}(t)$ satisfies the quadratic equation (in the variable $y$): 
$$ 4ty^2-(1+3t)y+1=0. 
\eqno(1)$$
The generating function $P_{\mathrm{rook}}(t)$ was derived earlier by  Coker \cite{Coker} by showing that it satisfies equation~(1).  In addition, Woan \cite{Woan} used a three-term recurrence to find $p_{\mathrm{rook},n}.$     

Little seems to have been done on Catalan queen paths.  For these paths, 
\begin{align*}  
T(t)  &= \frac {t}{1 - t}, 
\\   
C &= 1 - \frac {t}{1-t} \,=\, \frac {1-2t}{1-t},      
\\  
\alpha &=\frac{1+t-4t^2+(1-t)\sqrt{1-12t+16t^2}}{2t(2-3t)}. 
\end{align*}
 
\vskip 0.2in\noindent  
{\bf Theorem 3.3} (Catalan queen paths).  
$$
P_{\mathrm{queen}}(t)  = \frac{(1-t)(1+t-4t^2)-(1-t)^2\sqrt{1-12t+16t^2}}{2t(2-3t)^2},   
$$

\vskip 0.2in\noindent  
Expanding $P_{\mathrm{queen}}(t),$ we have 
$$
1+2t+10t^2+63t^3+454t^4+ 3539t^5 + 29008 t^6 + 246255 t^7 + 2145722 t^8 +  \cdots
$$
The first few coefficients appear as sequence A175962 in the OEIS, the \emph{On-line Encyclopedia of Integer Sequences}.  
We also have the following asymptotic formula.  

\vskip 0.2in\noindent  
{\bf Corollary 3.4} 
$$ p_{\mathrm{queen},n} \sim \frac{(35-15\sqrt{5})(3\sqrt{5}-5)^{1/2}}{2\sqrt{2}} \cdot \frac{(2(3+\sqrt{5}))^n}{\sqrt{\pi n^3}}. $$

\vskip 0.2in\noindent  
Note that $2(3+\sqrt{5})\approx 10.47$, while the growth constant of Catalan rook paths is $9$ and that of classical Catalan paths is $4$. 
From the formula for $P_{\mathrm{queen}}(t)$, we deduce that it satisfies the quadratic 
equation:
$$ t(2-3t)^2y^2-(4t^3-5t^2+1)y+(1-t)^2=0.$$

We next consider Catalan $I$-queen paths where $I$ is an initial segment of the positive integers.  
For a positive integer $r,$ let $S_r=\{(a,0),(0,a): 1 \leq a < \infty\} \cup \{(a,a): 1 \leq a \leq
r\}$ and $P_r(t)$ be the generating function for Catalan
$S_r$-paths.   Our method yields the following formulas for small values of $r:$   
\begin{align*}
P_1(t)&=\frac{1+2t-t^2-\sqrt{(t-1)(-1+11t-7t^2+t^3)}}{2t(t-2)^2}=1+2t+9t^2+57t^3+411t^4+\cdots \,\,,\\
P_2(t)&=\frac{1+3t+t^2-\sqrt{1-10t-5t^2+2t^3+t^4}}{2t(1-t)(2+t)^2}=1+2t+10t^2+62t^3+448t^4+\cdots \,\,,\\
P_3(t)&=\frac{1+2t-2t^2-2t^3-t^4-\sqrt{(t-1)(-1+11t-5t^2-5t^3-7t^4+t^5+t^6+t^7)}}{2t(-2+t+t^2+t^3)^2}\\
& =1+2t+10t^2+63t^3+453t^4+\cdots
\end{align*}
The first $r+1$-st terms of $P_r(t)$ and $P_{\mathrm{queen}}(t)$ agree, as expected.     The first few coefficients 
of $P_1(t)$ and $P_2(t)$ appear as sequences A175912 and A175939 in the OEIS.

A natural question is to ask for the generating functions for Catalan paths in which both rook and bishop steps have
bounded length.  As explained at the end of Section 2, the only cases that follow readily from our theory are those where $(1,0)$ is the only possible horizontal step.  
For a subset $J$ in $\mathbb{P},$ let $P_{1,J}(t)$ be the generating function for the number of Catalan $[\{(1,0),(0,1)\} \cup \{(a,a): a\in J\}]$-paths.
Using the same method as for Catalan $I$-queen paths, we obtain the following theorem. 

\vskip 0.2in\noindent  
{\bf Theorem 2.5.}  Let $J \subseteq \mathbb{P}.$  Then  
$$P_{1,J}(t)=\frac{1-T_J(t)-\sqrt{[1-T_J(t)]^2-4t}}{2t},$$
where $T_J(t) =\sum_{a\in J} t^a$.

\vskip 0.2in\noindent  
When $J = \emptyset$ and $J = \{1\},$ we have the classical generating functions for Catalan and Schr\"oder paths.  
%
%
For $J=\{2\},$ $J = \{1,2\},$ and $J=\mathbb{P} = \{1,2,\ldots\},$ we have
\begin{align*}
P_{1,\{2\}}(t)&=\frac{1-t^2-\sqrt{1-4t-2t^2+t^4}}{2t}=1+t+3t^2+ 8t^3+ 25t^4+\cdots
\\
P_{1,\{1,2\}}(t)&=\frac{1-t-t^2-\sqrt{1-6t-t^2+2t^3+t^4}}{2t}=1+2t+7t^2+27t^3+116t^4+\cdots
\\
P_{1,\mathbb{P}}(t)&=\frac{2t-1+\sqrt{1-8t+12t^2-4t^3}}{2t(t-1)}=1+2t+7t^2+28t^3+122t^4+\cdots
\end{align*}
The first few coefficients of $P_{1,\{1,2\}}(t)$ appear as sequence  A175934 in the OEIS.

\vskip 0.4in \noindent  
{\bf 4.  Catalan paths with queen and spider steps}
\medskip

For more general sets $S$ of steps, the conditions necessary to perform calculations akin to those in Section 3 are summarized in the following theorem. 

\medskip 
{\bf Theorem 4.1.\ }
Assume $S$ satisfies the slope condition for the boundary $(i+1)$ and that it
contains either all proper horizontal steps $(a,0)$, $a>0$, or exactly one 
horizontal step $(1,0).$ If the generating function $A(x,y)$ is rational and 
$T(t)$ is algebraic, then the generating function for Catalan $S$-paths is algebraic.   
%

\medskip

We will give two examples of paths to which the theorem applies.  
We begin with a simple example.  A {\sl Tugger} path is an $S$-path, where $S = \{(1,0),(0,1),(1,2)\}.$   The set $S$ contains the unit horizontal and vertical step, as well as the shortest spider step $(1,2).$  
It satisfies the slope condition for the Catalan boundary.  We can use Theorem 2.4 and the method in Section 3 to obtain the following formula: 
$$
P_{\mathrm{Tugger}}(t) = \frac {1-\sqrt{1-4t-4t^2}}{2t(1+t)}.
$$
Expanding $P_{\mathrm{Tugger}}(t),$ we have 
$$
1 + t + 3t^2 + 9t^3 + 31 t^4 + 113t^5 + 431 t^6 + 1697 t^7 + 6847 t^8 + \cdots.  
$$
Also, 
$$
p_{\mathrm{Tugger},n}\sim \sqrt{4-2\sqrt{2}}  \cdot \frac{(2+2\sqrt{2})^n}{\sqrt{\pi n^3}}. 
$$

Now let $M = \{(a,0),(0,a),(b,c):a   > 0, \, 1\leq b \leq c\}.$  As well as rook and queen steps, the set $M$ contains all spider steps and is the largest set of steps that satisfies the slope condition for the Catalan boundary.    For these paths, the bivariate generating function is given by 
$$
A_M(x,y)  = \left(2 - \frac {x}{1-x} - \left(\frac {1}{1-xy}\right)\left(\frac {1}{1-y}\right) \right)^{-1}.  
$$
Using Theorem 2.3 and the method in Section 3, we obtain the following formula for the generating function $P_M(t)$ for Catalan $M$-paths.

\vskip 0.2in\noindent 
{\bf Theorem 4.1} (Catalan $M$-paths).  
$$
P_M(t) =   \frac{ 3t^2 - t - 1 + \sqrt{1 - 14t +35t^2 - 30t^3 + 9t^4}}{4t(3t-2)} 
$$

\vskip 0.2in\noindent 
Expanding, we have 
$$
P(t) = 1 + 2t + 11t^2 + 75t^3 + 578t^4 + 4791t^5 + 41657 t^6 + 374728 t^7 + 3458073 t^8 + \cdots.  
$$
From the formula, we deduce that $P_M(t)$ satisfies the quadratic equation 
$$
2t(3t-2) y^2 - (3t^2 -t - 1)y + (t-1) = 0.  
$$ 
We also derive the asymptotic formula 
$$
p_{M,n} \sim \omega \frac{\gamma^n}{\sqrt{\pi  n^3}},
$$ 
where $\gamma$ and $\omega$ are algebraic numbers with $\gamma=11.0785\ldots$ and $\omega=0.6968\ldots$.

\vskip 0.4in \noindent  
{\bf 5. Step enumerators}
\medskip

Our methods can be used to calculate step enumerators of Catalan paths.  Let $S$ be a set of proper steps.  A {\sl weight} $w$  on $S$ is a function $S \to \mathbb{A},$ where $\mathbb{A}$ is a commutative ring with identity.  Usually, we will take $\mathbb{A}$ to be the ring of polynomials with complex coefficients in many variables and assign a variable as the weight of a step.  If $X$ is a set of $S$-paths, we define the {\sl step enumerator} $e^w[X]$ by 
$$
e^{w}[X] = \sum_{(s_1,s_2,\ldots,s_j) \in X}  w(s_1)w(s_2) \cdots w(s_j). 
$$   
The following lemma (which is almost tautological) says that step enumerators are {\sl multiplicative}.    

\vskip 0.2in\noindent 
{\bf Lemma 5.1.}  Let $X$ be a set of paths which can be decomposed into the cartesian product $Y \times Z$ so that each path in $X$ 
is the concatenation of a path in $Y$ and a path in $Z.$  Then 
$$
e^w[X] = e^w[Y] e^w[Z].
$$
 
\medskip\noindent 
Next, define $A^w(x,y)$ and $a_{n,m}^w$ by the equation:   
$$
A^w (x,y) = \sum_{n,m=0}^{\infty}  a_{n,m}^w x^ny^m
= \left( 1 - \sum_{(a,b) \in S}   w((a,b)) x^a y^b  \right)^{-1}.  
$$
Then $a_{n,m}^w$ is the step enumerator of the set of all $S$-paths from $(0,0)$ to $(n,m)$ and $A^w(x,y)$ is the bivariate generating function of $a_{n,m}^w.$  
In the same way as in Section 2, we define the generating functions $D_h^w(t)$ and $T^w(t).$    
Finally, let  $p^{w}_n$ be the step enumerator for Catalan $S$-paths ending at $(n,n)$ and    
$$
P^{w}(t) = \sum_{i=0}^{\infty}  p^{w}_n  t^n.  
$$

The proofs of the two theorems in Section~2 can easily be modified using Lemma 5.1 to yield the following weighted versions.  

\vskip 0.2in\noindent  
{\bf Theorem 5.2.}  

(a) If $S$ contains all proper horizontal steps $(a,0), a > 0,$ $S$ satisfies the slope condition for the boundary $(i+1)$, and the same weight $\rho$ is assigned to 
all horizontal steps $(a,0)$ in $S$ (regardless of length), then 
$$
P^{w} (t)=\frac{D_0^w(t)+D_1^w(t)+D_2^w(t)+\cdots}{1+(1 + \rho -T^w(t))(D^w_1(t)+D^w_2(t)+\cdots)}.   
$$

(b) If $S$ contains no horizontal step but $(1,0)$ and it is assigned the weight $\rho,$ then 
$$
P^{w} (t)=\frac{D_0^w(t)}{1 + \rho D_1^w(t)}.    
$$

\vskip 0.2in
We will illustrate Theorem 5.2 with the simple case of Tugger paths, where we impose the weights $w((1,0)) = w((0,1)) = \rho$ and $w((1,2)) = \sigma.$  Then 
$$
P^w_{\mathrm{Tugger}}(t;\rho,\sigma) = \frac {1 - \sqrt{1 - 4 \rho^2 t - 4 \rho \sigma t^2}}{2 \rho t (\rho + \sigma t)}.
$$
Expanding $P^w_{\mathrm{Tugger}}(t;\rho,\sigma) ,$ we have 
$$
1 + \rho^2 t + (2\rho^4 + \rho\sigma)t^2 + (5\rho^6 + 4 \rho^3\sigma)t^3 + (14 \rho^8 + 15\rho^5 \sigma + 2\rho^2 \sigma^2) t^4 + 
(42\rho^{10} + 56 \rho^7 \sigma + 15 \rho^4 \sigma^2) t^5 + \cdots .
$$

We next calculate generating functions for the step enumerators for Catalan rook, queen, and $M$-paths.  
These generating functions are rather complicated. The reader who needs to use them should calculate them using computer algebra.    For this purpose, we provide quadratic equations for which they are power-series roots.      
We begin with enumerators for rook paths by the number of horizontal and vertical steps. Let $w((a,0)) = \rho$ and $w((0,a)) = \nu.$    
Then 
\begin{align*}
P^{w}_{\mathrm{rook}}(t;\rho,\nu) &= \beta^w/(1 + \rho)  
\\
&= \frac {1 + (1 + \rho + \nu)t - \sqrt{1 - [1 + \rho + \nu + 2\rho\nu]t + [1+(\rho +\nu)(2+ \rho+\nu)]t^2}}   {2(1 + \rho)(1 + \nu)t} 
\end{align*}
where $\beta^w$  is the power-series root of the quadratic equation (in the variable $s$): 
$$
(1 + \nu)t s^2 - [1 + (1 + \rho + \nu)t ]s + (1+ \rho) = 0.  
$$
Expanding $P^{w}_{\mathrm{rook}}(t;\rho,\nu),$ we obtain  
$$
1 + \rho\nu t + (\rho\nu +\rho\nu^2 + \rho^2\nu + 2\rho^2\nu^2)t^2 + 
(\rho\nu +  2\rho\nu^2 + 2\rho^2\nu + \rho\nu^3 +\rho^3\nu + 7\rho^2\nu^2 + 5\rho^2\nu^3 + 5\rho^3\nu^2 + 5\rho^3\nu^3) t^3 + \cdots.  
$$
Setting $\nu = \rho,$  we obtain the generating function for enumerators of rook paths by the total number of steps:   
$$
P^w_{\mathrm{rook}}(t;\rho) =  \frac {1 + (1+2\rho)t - \sqrt{(1-t)(1- (1 + 2\rho)^2t)}}{2t(1+\rho)^2}. 
$$
 Expanding $P^{w}_{\mathrm{rook}}(t;\rho)$, we obtain  
$$
1 + \rho^2 t + (\rho^2 + 2 \rho^3 + 2 \rho^4) t^2 + (\rho^2 + 4\rho^3 + 9\rho^4 + 10\rho^5 + 5\rho^6 )t^3 + \cdots 
$$ 
Rook-path enumerators by total number of steps have been studied earlier in Coker \cite{Coker}.   

Moving on to queen paths, let a rook step have weight $\rho$ and a bishop step have weight $\omega.$  Then $P^w_{\mathrm{queen}}(t;\rho,\omega)$ equals 

%
$$
\frac {(1-t)[(1 - (\omega - 2\rho)t - (\omega + 2\rho + 1)t^2 ]
-(1-t)^2\sqrt{ 1 -[(2\rho+1)^2 + 2\omega + 1]t  + [1+ (2\rho + \omega + 2)(2\rho + \omega)]t^2 }}{2t[\rho + 1 -(\rho +\omega +1)t]^2}.
$$ 



\noindent
The generating function $P^w_{\mathrm{queen}}(t;\rho,\omega)$ is the power-series root of the quadratic equation (in the variable $y$):
$$
t[\rho + 1 -(\rho +\omega +1)t]^2 y^2 -(1-t)[(1 - (\omega - 2\rho)t - (\omega + 2\rho + 1)t^2 ]y + (1-t)^2 = 0.
$$
Expanding $P^w_{\mathrm{queen}}(t;\rho,\omega),$ we have 
$$
1 + (\rho^2 + \omega)t + (\omega +\omega^2 + 3\omega \rho^2 + \rho^2 + 2 \rho^3 + 2 \rho^4)t^2+ \cdots .
$$

We end with $M$-paths.   The actual generating functions are complicated.  In the simplest case, when all steps have the same weight $\rho,$  $P^w_M(t;\rho)$ is the power-series root of the quadratic equation:
$$
(1+\rho)t[(2\rho+1) t -\rho - 1] y^2 - [(2\rho+ 1) t^2 - \rho t - 1]y + (t-1) = 0.
$$
Expanding $P^w_M(t;\rho),$ we have 
$$
1 + (\rho + \rho^2)t + (\rho +3 \rho^2 + 5\rho^3 +2 \rho^4)t^2 + (\rho +5 \rho^2 + 17\rho^3 +27 \rho^4 + 20\rho^5 + 5 \rho^6)t^3 + \cdots .
$$


\vskip 0.2in \noindent  
{\bf 6. Other boundaries}
\medskip

Can our method for counting paths with Catalan boundaries be adapted to paths with other periodic boundaries?  To partially answer this question, consider $S$-paths with boundary $(2i+1)$, where $S$ is any set of steps that contains all proper horizontal steps and satisfies the slope condition for this boundary. From the recurrence in Corollary~2.2, we arrive at the following
expression for the generating function for these paths:
$$
P(t)= 
\frac {\displaystyle{\sum_{h=0}^{\infty} D_h(t)+ \sum_{h=1}^\infty \overline{D}_{-h}(t)}}
{\displaystyle{1+(2-T(t))\left[\sum_{h= 0}^\infty  D_h(t)+ \sum_{h =1}^\infty  \overline{D}_{-h}(t) -\sum_{n=0}^\infty a_{n,2n}t^n\right]}},
$$
where 
$$
\overline{D}_{-h}(t)=\sum_{n= h}^\infty  a_{n+h,n}t^n \quad \mathrm{and} \quad T(t)=\sum_{(2a,a)\in S} t^a. 
$$ 
Unfortunately, it seems hard to obtain a closed form for the series
$\overline{D}_{-h}(t)$ from the bivariate generating function $A(x,y)$.

\vskip 0.2in \noindent 
{\bf Acknowledgments.}  We were led to think about rook and queen paths with boundaries by the paper \cite{Kauers}, but we took an easier 
challenge.   JK was supported by the National Security Agency under Grant H98230-11-1-0183. AdM was
supported by the Spanish and Catalan governments  under projects  MTM2011-24097 and DGR2009-SGR1040.    
The calculations were done with the help of MAPLE.

\end{document}